\documentclass[12pt]{article}

\usepackage{amsmath,amsfonts,graphicx,color,fancyhdr,amssymb}

\title{Several Convex-Ear Decompositions}

\author{Jay Schweig \footnote{schweig@math.cornell.edu}}

\newtheorem{thm}{Theorem}[section]
\newtheorem{prop}[thm]{Proposition}
\newtheorem{cor}[thm]{Corollary}
\newtheorem{lem}[thm]{Lemma}
\newtheorem{defn}[thm]{Definition}
\newtheorem{fact}[thm]{Fact}
\newtheorem{conj}[thm]{Conjecture}

\begin{document} 

\maketitle

\begin{abstract}  In this paper we give convex-ear decompositions for the order complexes of several classes of posets, namely supersolvable lattices with non-zero M\"obius functions and rank-selected subposets of such lattices, rank-selected geometric lattices, and rank-selected face posets of shellable complexes which do not include the top rank.  These decompositions give us many new inequalities for the h-vectors of these complexes.  In addition, our decomposition of rank-selected face posets of shellable complexes allows us to prove inequalities for the flag h-vector of face posets of Cohen-Macaulay complexes. 
\end{abstract}

\section{Introduction}  

For a $(d-1)$-dimensional simplicial complex, its f-vector, $\langle f_0,f_1,\ldots ,f_d\rangle$, is the integral sequence that expresses the number of faces in each dimension.  This object is the primary combinatorial invariant associated to a simplicial complex $\Delta$.  Closely related to the f-vector is the h-vector $\langle h_0,h_1,\ldots, h_d\rangle$, defined by the equation $\Sigma_0^d f_i(x-1)^{d-i}=\Sigma_0^dh_ix^{d-i}$.  In the case when the complex to which it is associated is shellable, the h-vector has an interesting alternative, but still equivalent, interpretation.  A convex-ear decomposition, introduced by Chari in \cite{ch}, is an invaluable tool in showing several key inequalities of the h-vector.  

When a $(d-1)$-dimensional complex admits a convex-ear decomposition, its h-vector satisfies $h_i\leq h_{d-i}$ and $h_i\leq h_{i+1}$ for all $i<d/2$.  The g-vector of a complex is defined to be $\langle h_0, h_1-h_0, h_2-h_1, \ldots,  h_{\lfloor d/2 \rfloor}-h_{\lfloor d/2 \rfloor-1}\rangle$.  Following the celebrated g-theorem for simplicial polytopes (see \cite{g}), which states that the g-vector for a simplicial polytope is an M-vector (namely a degree sequence of some monomial order ideal), Swartz proved in \cite{sw} that the g-vector of any complex admitting a convex-ear decomposition is an M-vector.  This result is what one might expect, since complexes that admit convex-ear decompositions are essentially composed of polytopes attached to each other in a geometrically reasonable fashion.

In this paper we examine a few classes of poset order complexes that admit convex-ear decompositions.  In the order proven, these posets are supersolvable lattices with non-zero M\"obius function, rank-selected Boolean lattices, rank-selected supersolvable lattices with non-zero M\"obius function, rank-selected face posets of shellable complexes excluding the top rank level, and rank-selected geometric lattices.  Even though the first two posets are special cases of the third, we prove these separately since the proofs are in order of increasing complexity.  Also, because the methods of proof in all cases are similar, understanding one decomposition will help to understand the subsequent ones.  Finally, as a consequence of our convex-ear decomposition of rank-selected face posets of shellable complexes, we obtain inequalities for the flag h-vectors of face posets of Cohen-Macaulay complexes.   

Since all the posets mentioned above admit convex-ear decompositions, our main results concerning the h-vector can be summarized as follows.

\begin{thm}  Let $P$ be a poset that is either a supersolvable lattice with non-zero M\"obius function, a rank-selected subposet of such a lattice, a rank-selected subposet of a geometric lattice, or a rank-selected subposet of the face poset of some Cohen-Macaulay complex $\Sigma$ (where we require that the rank-selection does not include the elements of the face poset corresponding to facets of $\Sigma$).  Let $\Delta=\Delta(P-\{\hat{0},\hat{1}\})$ be the order complex of the proper part of $P$, and let $\langle h_0,h_1,\ldots,h_r\rangle$ be the h-vector of $\Delta$.  Then for $i<r/2$, 

1)  $h_i\leq h_{r-i}$, and

2)  $h_i\leq h_{i+1}$.

Furthermore, the g-vector of $\Delta$, $\langle h_0, h_1-h_0, h_2-h_1, \ldots, h_{\lfloor r/2\rfloor}-h_{\lfloor r/2\rfloor-1}\rangle$, is an M-vector.
\end{thm}

\section{Preliminaries}
\subsection{Shellability and Convex-Ear Decompositions}

Throughout this section, let $\Delta$ be a $(d-1)$-dimensional simplicial complex.  

For $0\leq i\leq d$, let $f_i$ be the number of simplices in $\Delta$ with $i$ vertices (so $f_0=1$ whenever $\Delta\neq\emptyset$).  Equivalently, $f_i$ is the number of $(i-1)$-dimensional simplices.  The \textit{h-vector} of $\Delta$ is the vector $\langle h_0,h_1,\ldots, h_d\rangle$ satisfying $\sum_0^df_i(x-1)^{d-i}=\sum_0^dh_ix^{d-i}$.  An \textit{M-vector} is the degree sequence of a monomial order ideal.  This is a natural generalization of the f-vector, since a simplicial complex can be viewed as a squarefree order ideal of monomials.

Two identities are immediate from this definition:  Plugging in $0$ for $x$, we get $h_d = \sum_0^df_i(-1)^{d-i}=(-1)^{d+1} \tilde{\chi}(\Delta)$.  Plugging in $1$ for $x$, we get $f_d=\sum_0^dh_i$.

\begin{defn}  $\Delta$ is \textit{shellable} (in the only sense considered in this paper) if it is pure and there exists an ordering $F_1, F_2,\ldots,F_t$ of its facets (i.e., maximal faces) such that for all $i$ with $1<i\leq t$, $F_i\cap (\bigcup_1^{i-1}F_j)$ is a non-empty union of facets of $\partial F_i$.
\end{defn}

The following is often given as an alternate definition of shellability (see for instance \cite{bw}):

\begin{prop}  \label{altshell}  An ordering $F_1,\ldots, F_t$ of the facets of $\Delta$ is a shelling if and only if for each $i,j$ with $i<j$, there exists a $k<j$ such that $F_i\cap F_j\subseteq F_k\cap F_j$ and $|F_k\setminus F_j|=1$.
\end{prop}

Suppose $F_1,F_2,\ldots,F_t$ is a shelling of $\Delta$.  For each $i>1$ let $r(F_i)=\{x\in F_i:F_i-x\subseteq \bigcup_1^{i-1}F_j\}$, and let $r(F_1)=\emptyset$.  Let $\Delta_k$ denote the simplicial complex generated by the facets $F_1, F_2,\ldots, F_k$, and note that each $r(F_i)$ has the following property:  if $K$ is a face of the simplicial complex $\Delta_i$ but not of the complex $\Delta_{i-1}$, then $r(F_i)\subseteq K$.  In other words, $r(F_i)$ is the unique minimal (with respect to inclusion) face in $\Delta_i\setminus \Delta_{i-1}$.

In the case when $\Delta$ is shellable, its h-vector can be interpreted as follows:

\begin{thm}
Suppose $\Delta$ has a shelling order $F_1,F_2,\ldots, F_t$ and let $h(\Delta)=\langle h_0,h_1,\ldots,h_d\rangle$.  Then for each $i$, $h_i=|\{F_j:|r(F_j)|=i\}|$.
\end{thm}   

The following proposition, proven by Danaraj and Klee, will be of use:

\begin{prop}  \label{ball}\cite{dk} Let $K$ be a full-dimensional shellable proper subcomplex of a sphere.  Then $K$ is a ball.  
\end{prop}

We are now in a position to define the concept of a convex-ear decomposition, originally introduced by Chari in \cite{ch}:

\begin{defn}  We say that a (d-1)-dimensional simplicial complex $\Delta$ has a \textit{convex-ear decomposition} if there exist subcomplexes
$\Sigma_1,\ldots\Sigma_n$ such that:

(i)  $\bigcup_1^n \Sigma_i = \Delta$.

(ii)  $\Sigma_1$ is the boundary complex of a simplicial $d$-polytope, while for
$i>1$ $\Sigma_i$ is a full-dimensional proper subcomplex of some simplicial
$d$-polytope.

(iii)  For $i>1$, $\Sigma_i$ is a simplicial ball.

(iv)  For $i>1$, $(\bigcup_1^{i-1}\Sigma_j)\cap \Sigma_i=\partial \Sigma_i$.
\end{defn}

\begin{thm} \cite{ch} Let $\Delta$ be a $d-1$-dimensional simplicial complex that
admits a convex-ear decomposition.  Then for $i<d/2$ the $h$-vector of $\Delta$
satisfies:  

1) $h_i\leq h_{d-i}$, and 

2) $h_i\leq h_{i+1}$.
\end{thm}

Swartz has also proved (in \cite{sw}) the following analogue of Stanley's g-theorem:

\begin{thm}  Let $\Delta$ be as in the statement of the previous theorem.  Then the
$g$-vector of $\Delta$, $\langle h_0, h_1-h_0, h_2-h_1, \ldots, h_{\lfloor d/2
\rfloor}-h_{\lfloor d/2 \rfloor-1}\rangle$, is an $M$-vector.
\end{thm}

Let $\Delta$ be a simplicial complex, and for any face $F$ of $\Delta$ define the link of $F$, written $lk(F)$, to be the subcomplex $\{G\in\Delta: F\cap G=\emptyset$ and $F\cup G\in \Delta\}$.  We say that $\Delta$ is \textit{Cohen-Macaulay} if, for every face $F$ of $\Delta$, $\tilde{H}_i(lk(F))=0$ whenever $i<\dim lk(F)$.  We leave it to the reader to verify the following implication:

\begin{prop}  If $\Delta$ admits a convex-ear decomposition or is shellable, then $\Delta$ is Cohen-Macaulay.
\end{prop}

For a vertex $v$ of $\Delta$, define the simplicial complex $\Delta-v$ to be the simplicial complex with faces of the form $F-v$ whenever $F\in \Delta$.  We call a Cohen-Macaulay complex $\Delta$ \textit{2-Cohen Macaulay}, or \textit{2-CM}, if for every vertex $v$ of $\Delta$ $\dim(\Delta-v)=\dim(\Delta)$ and $\Delta-v$ is Cohen-Macaulay.

\begin{prop}  If $\Delta$ admits a convex-ear decomposition, then $\Delta$ is 2-CM.
\end{prop}    

\subsection{Order Complexes}
Proofs of the following can be found in \cite{bw}.

Let $P$ be any finite poset.  The \textit{order complex} of $P$, written $\Delta(P)$, is the simplicial complex whose simplices are chains in $P$.  Note that $\Delta(P)$ is pure if and only if $P$ is graded (i.e., ranked).

Associated to any simplicial complex $K$ is its \textit{face poset} $P_K$.  The elements of $P_K$ are the faces of $K$, and the order relation is inclusion.  The proof of the following is left as an easy exercise for the reader:

\begin{prop}  Let $K$ be a simplicial complex.  Then the order complex of its face poset, $\Delta(P_K-\emptyset)$, is the first barycentric subdivision of $K$.
\end{prop}   

Thus a simplicial complex and the order complex of its face poset are homeomorphic.

Let $P$ be a graded poset with greatest element $\hat{1}$ and least element $\hat{0}$, and let $\epsilon(P)=\{(x,y)\in P^2: y$ covers $x\}$ be the set of edges of the Hasse diagram of $P$.  An \textit{EL-labeling} of $P$ is a function $\lambda:\epsilon(P)\rightarrow \mathbb{N}$ satisfying the following:

i)  In each interval $[x,y]$, there is a unique saturated chain $\textbf{c}:=x=x_0<x_1<\ldots<x_r=y$ such that the $r$-tuple $\lambda(\textbf{c})$, by which we mean $\langle \lambda(x_0,x_1), \lambda(x_1,x_2), \ldots, \lambda(x_{r-1},x_r)\rangle$, is weakly increasing.

ii)  If $\textbf{d}$ is any other saturated chain in the interval $[x,y]$, then $\lambda(\textbf{c})$ lexicographically precedes $\lambda(\textbf{d})$.  

\begin{thm}
If $P$ is a graded poset that admits an EL-labeling, then lexicographic order of the maximal chains in $P$ corresponds to a shelling of $\Delta(P-\{\hat{0},\hat{1}\})$.
\end{thm}   

Whenever we say $P$ admits an EL-labeling, we assume implicitly that $P$ has a greatest element $\hat{1}$ and a least element $\hat{0}$ and that $P$ is graded.

For any poset $P$, the M\"obius function $\mu:P^2\rightarrow \mathbb{Z}$ is defined recursively by $\mu(x,x)=1$, $\mu(x,y)=-\sum_{x\leq z<y}\mu(x,z)$ if $x<y$, and $\mu(x,y)=0$ if $x\nleq y$.  

\begin{prop}\label{mu}  (see for instance \cite{st3})  Let $P$ be a poset that admits some EL-labeling, and let $\lambda$ be any EL-labeling of $P$.  Then $|\mu(x,y)|$ is the number of maximal chains $\textbf{c}$ in the interval $[x,y]$ such that $\lambda(\textbf{c})$ is strictly decreasing.
\end{prop}

\begin{prop}\label{mu2}  Let $P$ be a poset with a least element $\hat{0}$ and a greatest element $\hat{1}$.  Then $\mu(\hat{0},\hat{1})=\tilde{\chi}(\Delta(P-\{\hat{0},\hat{1}\}))$.
\end{prop}

Let $\sigma=a_1a_2\ldots a_m$ be a word of integers.  The \textit{descent set} of $\sigma$ is $d(\sigma)=\{i:a_i>a_{i+1}\}$.  It should be noted that the following, along with the previous proposition, implies Proposition \ref{mu}:  

\begin{prop}\label{d}  Let $P$ be a poset that admits an EL-labeling $\lambda$, and let $\Delta$ be the order complex of the proper part of $P$.  Then $h_i(\Delta)=|\{$maximal chains $\textbf{c}$ of $\Delta: |d(\lambda(\textbf{c}))|=i\}|$.
\end{prop}

\subsection{Flag Vectors and the Weak Order}

Let $P$ be a graded poset of rank $r$, and let $\Delta$ be its order complex.  For any subset $S\subseteq 
[r]$, define 
$f_S$ to be the 
number of maximal chains of the rank selected subposet $P_S=\{x\in P:$ rank$(x)\in S\}$.  This gives a natural refinement of the f-vector of $\Delta$ since $f_i=\sum_{|S|=i}f_S$.

We also define the flag h-vector of $\Delta$ by $h_S=\sum_{T\subseteq S}(-1)^{|S-T|}f_T$.  Equivalently, by inclusion-exclusion (see \cite{st3}), $f_S=\sum_{T\subseteq S}h_T$.  If the poset $P$ admits an EL-labeling, the flag h-vector provides a nice enumerative interpretation:

\begin{prop}(see for instance \cite{st3}, pg. 133)  Let $P$ be a poset that admits an EL-labeling.  Then $h_S$ counts the number of maximal chains of $P$ whose labels have descent set $S$.
\end{prop}  

The above proposition in conjunction with Proposition \ref{d} implies that whenever $P$ admits an EL-labeling the flag h-vector of $P$ satisfies:

$$h_i=\sum_{|S|=i}h_S$$

In fact, this is true for any graded poset.  Thus the flag h-vector is a refinement of the usual h-vector in the same way that the flag f-vector is a refinement of the usual f-vector.

Now let $\sigma\in S_m$ be a permutation written as a word in $[m]$: $\sigma=\sigma_1 \sigma_2 \ldots \sigma_m$.  If we interchange $\sigma_i$ and $\sigma_{i+1}$ for some $i\notin d(\sigma)$, we call this a \textit{switch}, and we say that $\sigma$ is less than $\tau$ in the \textit{weak order} (sometimes called the \textit{weak Bruhat order}), written $\sigma<_w\tau$, if $\tau$ can be obtained from $\sigma$ by a sequence of switches. 

For example, $1\hspace{2pt}3\hspace{2pt}2\hspace{2pt}4 <_w 1\hspace{2pt} 4 \hspace{2pt} 2 \hspace{2pt} 3$, since $1 \hspace{2pt} 3 \hspace{2pt} 2 \hspace{2pt} 4<_w 1 \hspace{2pt} 3 \hspace{2pt} 4 \hspace{2pt} 2<_w 1 \hspace{2pt} 4 \hspace{2pt} 3 \hspace{2pt} 2$.

For a subset $S\subseteq[m]$, let $D_S=\{\sigma\in S_m: d(\sigma)=S\}$.  For two subsets $S,T\subseteq[m]$, we say that $S$ \textit{dominates} $T$ if there exists an injection $\phi:D_T\rightarrow D_S$ such that $\tau<_w\phi(\tau)$ for all $\tau\in D_T$.  For example, if $m=4$, the set $\{1,3\}$ dominates the set $\{1\}$, since for any permutation $\tau\in D_{\{1\}}$, $\phi(\tau)=\tau_1 \tau_2 \tau_4 \tau_3\in D_{\{1,3\}}$, and this map is clearly injective.  For a study of which pairs of subsets $S,T\subseteq [m]$ satisfy this dominance relation, see \cite{de} or \cite{ns}.    

\section{The Supersolvable Lattice}\label{super}

Supersolvable lattices originally arose from group theory as subgroup lattices of supersolvable groups.  However, results like that of Stanley's (\cite{st1}) and McNamara's (\cite{mc}) show that supersolvable lattices are also of interest from a purely combinatorial perspective.  Their topology has also been investigated by Welker, who has shown (\cite{we}) that the order complex of a supersolvable lattice with non-zero M\"obius function is 2-Cohen-Macaulay.  The following combinatorial definition was formulated by Stanley in \cite{st1}:

\begin{defn}  A lattice $L$ is supersolvable if there is some maximal chain $\textbf{c}$ (called the \textit{M-chain} of $L$; not to be confused with an M-vector) such that the sublattice of $L$ generated by $\textbf{c}$ and any other (not necessarily maximal) chain $\textbf{d}$ is distributive.
\end{defn}

\begin{thm}\label{sslattice}  Let $L$ be a supersolvable lattice such that whenever $x,y\in
L$ and $x<y$ then $\mu(x,y)\neq 0$, and let
$\Delta=\Delta(L\setminus\{\hat{0},\hat{1}\})$ be the order complex of the proper
part of $L$.  Then $\Delta$ admits a convex-ear decomposition.
\end{thm}

Now let $L$ be as in the hypothesis of the theorem, and suppose that $L$ has rank $r$.  Stanley
has shown (\cite{st1}) that $L$ admits an $S_r$-labeling, which is an EL-labeling with labels
from $[r]$ such that no two edges in the Hasse diagram of $L$ that are in some chain
together have the same label.  In this labeling, the M-chain of $L$ is the unique maximal chain with increasing labels.

Let $\lambda$ be an $S_r$-labeling of $L$, and let $[x,y]$ be an interval in $L$.  By
definition of an EL-labeling, there is a unique saturated chain
$x=a_1<a_2<\ldots<a_{{\rm rank}(y)-{\rm rank}(x)+1}=y$ such that
$\lambda(a_{i-1},a_i)<\lambda(a_i,a_{i+1})$ for all $i$.  Furthermore, since $\mu(x,y)\neq 0$, Proposition \ref{mu} gives us:

\begin{fact}  In every interval of $L$, there is a unique saturated chain with increasing labels,
and there is at least one saturated chain with decreasing labels.  
\end{fact}

Let $\textbf{c}$ be the M-chain of $L$.  That is, $\textbf{c}$ is the unique
saturated chain in the interval $ [ \hat{0} , \hat{1} ] $ (i.e., the whole lattice)
with increasing labels.  We know by the above fact that there is at least one
chain in $[\hat{0},\hat{1}]$ with decreasing labels.  Let $\textbf{c}_1 ,\ldots
,\textbf{c}_t$ be all such chains, and for $i=1,\ldots, t$, let $L_i$ be the
sublattice of $L$ generated by $\textbf{c}$ and $\textbf{c}_i$.  By definition, each
$L_i$ is a distributive lattice.  

We claim that each $L_i$ is isomorphic to $B_r$, the boolean lattice on $r$
elements.  To see this, let $P$ be the $r$-element poset for which $L_i$ is the
lattice of order ideals (such a poset is guaranteed by the fundamental theorem of finite distributive lattices; see for instance \cite{st3}).  Now $\textbf{c}$ gives us a chain of order ideals of $P$:  $\emptyset =
I_0\subset I_1\subset \ldots \subset I_r = P$.  Note that $|I_j\setminus
I_{j-1}|=1$, since $\textbf{c}$ is a saturated chain.  For $j=1,\ldots r$, let
$x_j=I_j\setminus I_{j-1}$.  Then $\textbf{c}$ corresponds to an order completion of
$P$, namely $x_1<x_2<\ldots x_r$, while $\textbf{c}_i$ corresponds to another order
completion of $P$:  $x_r<x_{r-1}<\ldots <x_1$.  Thus, $P$ must be the $r$-element
antichain, and so every subset of $P$ is an order ideal.  The isomorphism $L_i\simeq
B_r$ is now obvious.  It is also known (see for instance \cite{st1}), that the given
$S_r$-labeling when restricted to the sublattice $L_i$ is the same as the labeling
given by viewing $L_i$ as the lattice of order ideals of $P$.  We now know
that for each $L_i$ and any permutation $\pi:[r]\rightarrow [r]$, there is a unique
saturated chain $\hat{0}=x_0<x_1<\ldots<x_r =\hat{1}$ in $L_i$ such that
$\lambda(x_{j-1},x_j)=\pi(j)$ whenever $1\leq j\leq r$.  

For $1\leq i\leq t$, let $E_i$ consist of all maximal chains of $L_i$ that are not
maximal chains in any $L_j$ for any $j<i$, let
$\hat{E}_i=\{\textbf{d}\setminus\{\hat{0},\hat{1}\}:\textbf{d}\in E_i\}$, and let
$\Sigma_i$ be the simplicial complex whose facets are given by the chains of
$\hat{E}_i$.  

Given these new constructions, we restate Theorem \ref{sslattice} in more specific terms as a proposition:

\begin{prop}  $\Sigma_1,\Sigma_2,\ldots,\Sigma_t$ is a convex-ear decomposition
of $\Delta$. 
\end{prop}

First note that $E_1$ is just the set of all maximal chains in $L_1$, and so $\Sigma_1$ is
really $\Delta (B_r\setminus \{\hat{0},\hat{1}\})$, the order complex of the proper
part of the boolean lattice on $r$ elements.  Thus, $\Sigma_1$ is the first
barycentric subdivision of the boundary of the $(r-1)$-simplex, meaning it is the
boundary complex of a simplicial $(r-2)$ polytope.  Because $\textbf{c}\notin E_i$ for $1<i\leq t$, each $\Sigma_i$ is a
proper subcomplex of $\Delta(L_i\setminus \{\hat{0},\hat{1}\})=\Delta (B_r\setminus
\{\hat{0},\hat{1}\})$, which is the boundary complex of a simplicial
$(r-2)$-polytope.  This proves property (ii) of the decomposition.

Next we show property (i).  Let $\textbf{d}:=\hat{0}=x_0<x_1<\ldots<x_r =\hat{1}$ be
any saturated chain of $L$.  If $\lambda(x_{i-1},x_i)>\lambda(x_i,x_{i+1})$ for all
i, then $\textbf{d}=\textbf{c}_j$ for some $j$.  Otherwise, there is an $i$ such
that $\lambda(x_{i-1},x_i)<\lambda(x_i,x_{i+1})$.  Since $\mu(x_{i-1},x_{i+1})\neq
0$, there exists $y\in L$ with $x_{i-1}<y<x_{i+1}$ such that
$\lambda(x_{i-1},y)>\lambda(y,x_{i+1})$.  Let $\textbf{d}':=
\hat{0}=x_0<x_1<\ldots<x_{i-1}<y<x_{i+1}<\ldots <x_r =\hat{1}$, and assume that
$\textbf{d}'$ is a chain in $L_j$ for some $j$.  Let $\pi\in S_r$ be the
permutation defined by $\pi(k)=\lambda(x_{k-1},x_k)$.  Then $L_j$ must contain a
saturated chain with $\pi$ as its labels, meaning $L_j$ contains a chain of the form
$\hat{0}=x_0<x_1<\ldots<x_{i-1}<z<x_{i+1}<\ldots <x_r =\hat{1}$ satisfying
$\lambda(x_{i-1},z)=\lambda(x_{i-1},y)$ and $\lambda(z,x_{i+1})=\lambda(y,x_{i+1})$.
 Since $\lambda(x_{i-1},y)< \lambda(y,x_{i+1})$, uniqueness of the increasing
chain in $[x_{i-1},x_{i+1}]$ gives us that $z=y$, meaning $\textbf{d}$ is a chain in
$L_j$.  

So for any saturated chain $\textbf{d}$ with an ascent, we can create a chain
$\textbf{d}'$ such that $\textbf{d}'$ is lexicographically later than
$\textbf{d}$, and such that if $\textbf{d}'$ is a chain in some $L_j$ then
$\textbf{d}$ is a chain in $L_j$ as well.  Therefore, starting with a $\textbf{d}$, we can
apply this process repeatedly until we reach a chain with no ascents.  This chain
will then be $\textbf{c}_k$ for some $k$, meaning our original chain $\textbf{d}$ is
in the sublattice $L_k$.  Hence, $\textbf{d}$ is in $E_n$ for some
$n\leq k$, and the facet corresponding to $\textbf{d}$ in $\Delta$ is
contained in $\Sigma_n$.  So $\Delta=\bigcup_1^t \Sigma_i$, and we have proven
property (i) of the decomposition. 

We now wish to show that $\Sigma_i$ is shellable whenever $i>1$.  This, coupled with Proposition \ref{ball}, will imply property (iii) of the decomposition. 

\begin{lem}  Fix $i$, let $\textbf{d}_1:=\hat{0}=x_0<x_1<\ldots<x_r=\hat{1}$ and
$\textbf{d}_2:=\hat{0}=y_0<y_1<\ldots<y_r=\hat{1}$ be two maximal chains of $L_i$, and
suppose that $\textbf{d}_1$ is lexicographically earlier than $\textbf{d}_2$.  Then
there exists a $j\in\{1,2,\ldots, r-1\}$ such that $x_j\neq y_j$ and
$\lambda(x_{j-1},x_j)<\lambda(x_j,x_{j+1})$.
\end{lem}

\textbf{Proof:}  This is a direct consequence of the fact that $\lambda$ is an EL-labeling.  $\square$ 

We are now in a position to shell each $\Sigma_i$.  We do this by invoking Proposition \ref{altshell}.  Fix $i>1$, and
let $\textbf{d}_1,\ldots \textbf{d}_n$ be the maximal chains of $E_i$ listed in
reverse lexicographic order of their $\lambda$-labels.  Note that $\textbf{d}_1=\textbf{c}_i$, and suppose $\textbf{d}_\alpha:=\hat{0}=y_0<y_1<\ldots<y_r=\hat{1}$ is earlier in the order than
$\textbf{d}_\beta:=\hat{0}=x_0<x_1<\ldots<x_r=\hat{1}$.  Then $\textbf{d}_\beta$ is
lexicographically earlier than $\textbf{d}_\alpha$.  By the lemma, we can find a $j\in
[r]$ so that the conclusion of the lemma holds.  Now since $\mu(x_{j-1},x_{j+1})\neq
0$, we can find a $z\in (x_{j-1},x_{j+1})$ satisfying
$\lambda(x_{j-1},z)>\lambda(z,x_{j+1})$.  Let $\textbf{d}:=\hat{0}=x_0<x_1<\ldots
<x_{j-1}<z<x_{j+1}<\ldots<x_r=\hat{1}$.  We claim that $\textbf{d}\in E_i$.  Indeed,
suppose that $\textbf{d}\in E_k$ for some $k<i$.  Then $\textbf{d}$ would be a
maximal chain in $L_k$.  However, uniqueness of the increasing chain in the interval
$[x_{j-1},x_{j+1}]$ would then give us that $\textbf{d}_\beta$ is a chain in $L_k$, a
contradiction since $\textbf{d}_\beta\in E_i$.  So $\textbf{d}\in E_i$,
$|\textbf{d}_\beta\setminus \textbf{d}|=1$, and $\textbf{d}_\beta\cap \textbf{d}_\alpha \subseteq \textbf{d}_\beta \cap \textbf{d}$.  Since $\textbf{d}$ is lexicographically later (and therefore earlier in the order) than $\textbf{d}_\beta$, this ordering is a shelling
of $\Sigma_i$.  This proves part (iii) of the decomposition.  $\square$

Finally, we need to show property (iv) of the decomposition.  Fix $i>1$, and note
that a non-maximal chain $\textbf{d}$ represents a simplex in $\partial \Sigma_i$ if
and only if:

1) $\textbf{d}$ is a subchain of some maximal chain in $\hat{L}_i\setminus E_i$,
where $ \hat{L}_i$ is the set of all maximal chains in $L_i$ (i.e., $\textbf{d}$ is
a subchain of an old maximal chain), and: 

2) $\textbf{d}$ is a subchain of some maximal chain in $E_i$ (i.e., $\textbf{d}$ is
a subchain of a new maximal chain).

First, we note that 1) immediately gives $\partial \Sigma_i\subseteq
\Sigma_i\cap (\bigcup_1^{i-1}\Sigma_j)$.  Let $p<r$, and let $\textbf{d}:=x_0<x_1<\ldots<x_p$ be some non-maximal chain representing a simplex in $ \Sigma_i\cap
(\bigcup_1^{i-1}\Sigma_j)$.  Since $\textbf{d}$ represents a simplex in $\Sigma_i$, it must satisfy condition 2) above.  This simplex is also contained in $\bigcup_1^{i-1}\Sigma_j$,  meaning $\textbf{d}$ must be a chain in $L_j$ for some $j<i$.  Now take every `gap' in
$\textbf{d}$ (i.e., every interval $[x_k,x_{k+1}]$ such that rank$(x_k)\neq
$ rank$(x_{k+1})-1$) and fill it with the unique chain in that interval with increasing
labels.  By uniqueness of these chains, the resulting maximal chain must be in
$L_j$, since $\textbf{d}$ is.  Therefore, $\textbf{d}$ satisfies 1), and we have
shown the final part of the decomposition.  $\square$

\section{The Rank-Selected Boolean Lattice}\label{rankbool}

Let $B_r$ denote the Boolean lattice of rank $r$, and let $\lambda:\{(x,y)\in B_r^2:y$ covers $x\}\rightarrow [r]$ be the standard EL-labeling on $B_r$, where if we
view $x$ and $y$ as subsets of $[r]$ with $y=x\cup\{k\}$ then $\lambda(x,y)=k$. 
It is clear that every element $\sigma\in S_n$ shows up exactly once as the labeling of a
maximal chain in $B_n$, and that any two saturated chains beginning with $y\in B_r$ and ending with $z\in B_r$ are labeled with the same subset of $[r]$.  Call this last subset $\lambda(y,z)$.\\

For the remainder of this section, fix a subset $S\subseteq [n-1]$.  Let $D=\{\sigma\in S_r: d(\sigma)=S\}$.  Let $\sigma_1,\sigma_2,\ldots,
\sigma_t$ be the elements of $D$ in lexicographic order, and let
$\textbf{d}_{\sigma_1} :=\textbf{d}_1,\ldots,\textbf{d}_{\sigma_t}:=\textbf{d}_t$ be
the corresponding maximal chains.  Fix $\sigma_i\in D$, and let
$\textbf{d}_i:=\hat{0}=x_0<x_1<\ldots<x_r=\hat{1}$.  Define $C_i=
\{\textbf{c}:=\hat{0}=y_0<y_1<\ldots<y_{r-1}<y_r=\hat{1}:j\notin S\Rightarrow
y_j=x_j\}$.\\

For each $i$ with $1\leq i\leq t$, define the sublattice $L_{\sigma_i}:=L_i$ as the set $\{z\in L:z\in \textbf{c}$ for some $\textbf{c}\in C_i\}$ with order relations inherited from $L$, and let $E_i$ be the set of all maximal chains belonging to $L_i$ that are not maximal chains in $L_j$ for any $j<i$.  Also, let $\Sigma_i$ be the simplicial complex with facets given by the maximal chains in $E_i-\{\hat{0},\hat{1}\}$.

\begin{thm}  $\Sigma_1,\Sigma_2, \ldots, \Sigma_t$ is a convex-ear decomposition of
$\Delta(L)$.  
\end{thm}

Write $S=\beta_1\cup\ldots\cup\beta_\ell$ where each $\beta_j$ is an interval in the
integers, $\beta_j\cup\beta_k$ is not an interval in the integers for any choice of
$j$ and $k$, and $\min\beta_1<\min\beta_2<\ldots<\min\beta_\ell$.  If we let
$L_i^j=\{z\in L_i:$ rank$(z)\in\beta_j\}$, then $L_i^j$ is isomorphic to the Boolean
lattice $B_{|\beta_j|+1}$ minus the top and bottom elements, meaning that its order
complex $\Delta(L_i^j)$ is the first barycentric subdivision of
$\partial\Delta^{|\beta_j|}$.  Since any selection of one maximal chain from each of
$L_i^1,L_i^2,\ldots L_i^\ell$ gives a maximal chain in $L_i$,
each $\Sigma_i$ is a subcomplex of $\Delta(L_i^1)\ast\ldots\ast\Delta(L_i^\ell)$, where ``$\ast$" denotes simplicial join.  It follows that each $\Sigma_i$ triangulates the boundary of a $|S|-$polytope.  This verifies property (ii) of the decomposition.  $\square$

For a set $A$ of integers, let $u(A)$ denote the word that is the elements of $A$
written in increasing order, and let $v(A)$ denote the word that is the elements of
$A$ written in decreasing order (if $A=\emptyset$ then we let $u(A)=v(A)$ be the
empty word).  Also, for each $\beta_i$ (as defined previously) let
$a_i=\min\beta_i$ and $b_i=\max\beta_i$.  That is, $\beta_i=[a_i,b_i]$.  We ignore the case when $a_i=b_i$ for some $i$, as this is an easy generalization from the case we treat.  Now let $\textbf{c}:=x_{a_1}<x_{a_1+1}<\ldots<x_{b_1}<x_{a_2}<x_{b_2}<\ldots<x_{b_\ell}$ be
a maximal chain in $L$ (where the indices of the $x_i$'s denote their ranks), and let
$\textbf{c}':=\hat{0}=x_0<x_1<\ldots<x_n=\hat{1}$ be the maximal chain obtained from
$\textbf{c}$ by filling in each gap in $\textbf{c}$ with the unique maximal chain in
that interval with increasing labels.  For $i,j\in[r]$ with $i\leq j$, let
$\lambda(i,j)$ denote the set
$\{\lambda(x_i,x_{i+1}),\lambda(x_{i+1},x_{i+2}),\ldots,\lambda(x_{j-1},x_j)\}$.  If
$i=j$, set $\lambda(i,j)=\emptyset$.  

Define $\sigma^{\mathbf c}$ to be the permutation
$u(\lambda(0,a_1-1))v(\lambda(a_1-1,b_1+1))u(\lambda(b_1+1,a_2-1))\ldots
u(\lambda(b_\ell+1,r))$.  The reason for introducing $\sigma^{\mathbf c}$ is the following:

\begin{lem}  \label{lex}$\sigma^{\mathbf c}$ is lexicographically least in the set $\{\tau\in
D:\textbf{c}$ is a maximal chain in $L_\tau\}$.
\end{lem}

\textbf{Proof:}  In order for the proof to cover all cases, let $b_0=0$ and
$a_{\ell+1}=r$.  Any $\tau\in D$ such that $\textbf{c}$ is a maximal chain in
$L_\tau$ must satisfy the following property:  For each $i\in[\ell]$, the set
$\{\tau(a_i),\tau(a_i+1),\ldots,\tau(b_i+1)\}$ must consist of the set
$\lambda(a_i,b_i)$ as well as one element each from the sets $\lambda(b_{i-1},a_i)$
and $\lambda(b_i,a_{i+1})$.  Since in the case of $\sigma^{\mathbf c}$ these elements are
$\lambda(x_{a_i-1},x_{a_i})=\max\lambda(b_{i-1},a_i)$ and
$\lambda(x_{b_i},x_{b_i+1})=\min\lambda(b_i,a_{i+1})$, respectively, the proof of
the lemma follows.  $\square$

The existence of $\sigma^{\mathbf c}$ shows that every maximal chain in $L$ is in
$E_i$ for some $i\in[t]$.  This gives us property (i).

Now fix $i\in[t]$, let $\textbf{d}_i:=\hat{0}=y_0<y_1<\ldots<y_r=\hat{1}$, and let
$\textbf{c}:=x_{a_1}<x_{a_1+1}<\ldots<x_{b_1}<x_{a_2}<\ldots<x_{b_\ell}$ be a
maximal chain in $L_i$.  The label $\lambda(\textbf{c})$ is defined to be the
following word:

$\lambda(y_{a_1-1},x_{a_1})\lambda(x_{a_1},x_{a_1+1}) \ldots
\lambda(x_{b_1},y_{b_1+1}) \lambda(y_{a_2-1},x_{a_2}) \ldots
\lambda(x_{b_\ell},y_{b_\ell+1})$. 

The following lemma should remind the reader of our shelling when dealing with the
convex-ear decomposition of a supersolvable lattice's order complex:

\begin{lem} [property (iii)]  Reverse lexicographic order of the maximal chains in
$E_i$ is a shelling of $\Sigma_i$.
\end{lem}

\textbf{Proof:}  Let $\textbf{c}_i$ denote the restriction of $\textbf{d}_i$ to all elements with ranks in $S$.  It is clear by construction then that
$\sigma^{{\mathbf c}_i}=\sigma_i$.  Therefore, the lexicographically greatest maximal chain in $L_i$ is always in $E_i$.  Let $\textbf{c}_j$ and $\textbf{c}_k$ be two maximal
chains in $E_i$, with $\textbf{c}_j$ earlier in the ordering than $\textbf{c}_k$. Write $\textbf{c}_k$ as: 
$$y_{a_1-1}<x_{a_1}<x_{a_1+1}<\ldots<x_{b_1}<x_{b_1+1}<x_{a_2-1}<x_{a_2}<\ldots<x_{b_\ell}<y_{b_\ell+1}$$ 

Each element of the form $y_{a_i-1}$ or $y_{b_i+1}$ is not really in $\textbf{c}_k$, but we write them for heuristic reasons.  Now $\textbf{c}_j$ is lexicographically later than $\textbf{c}_k$, and so we are guaranteed the existence of some $m\in S$ such
that $x_m\notin \textbf{c}_j$ and $\lambda(x_{m-1},x_m)<\lambda(x_m,x_{m+1})$.  Since the open interval $(x_{m-1},x_{m+1})$ has cardinality $2$, let $x_m'$ be the
element in the interval other than $x_m$.  Then
$\lambda(x_{m-1},x_m')=\lambda(x_m,x_{m+1})$ and
$\lambda(x_m',x_{m+1})=\lambda(x_{m-1},x_m)$.  Let $\textbf{c}_k'$ be the chain that
results from replacing $x_m$ with $x_m'$ in $\textbf{c}_k$.  If $a_n<m<b_n$ for
some $n\in[\ell]$, it is immediate from the construction of $\sigma^{{\mathbf c}_k}$ that
$\sigma^{{\mathbf c}_k'}=\sigma^{{\mathbf c}_k}=\sigma_i$, and so $\textbf{c}_k'\in E_i$.  Suppose that $m=b_n$ for some $n\in[\ell]$.  Because $\sigma^{{\mathbf c}_k}=\sigma_i$, $\lambda(x_m,y_{b_n+1})$ is less then each element of the set
$\{\lambda(y_{b_n+1},y_{b_n+2}),\lambda(y_{b_n+2},y_{b_n+3}),\ldots,
\lambda(y_{a_{n+1}-1},x_{a_{n+1}})\}$.  On the other hand,
$\lambda(x_m',y_{b_n+1})=\lambda(x_{m-1},x_m)< \lambda(x_m,x_{m+1})$, so $\lambda(x_m',y_{b_n+1})$ is also less than every element in the above set. 
The construction of $\sigma^{{\mathbf c}_k'}$ then mirrors the construction of $\sigma^{{\mathbf c}_k}$,
meaning $\sigma^{{\mathbf c}_k'}=\sigma_i$ and $\sigma^{{\mathbf c}_k'}\in E_i$.  The case when $m=a_n$
is completely symmetric.  $\square$

Implicit in the previous proof is the following fact:

Fix $i$, let $\textbf{c}':=\hat{0}=x_0<x_1<\ldots<x_r=\hat{1}$ be a maximal chain in
$C_i$, and let $\textbf{c}$ be its restriction to elements with ranks in $S$.  It is clear that $\textbf{c}'$ is completely determined by $\textbf{c}$.  Furthermore, let
$$\lambda(n,s)=\{\lambda(x_n, x_{n+1}), \lambda(x_{n+1}, x_{n+2}),  \ldots,  \lambda(x_{s-1},x_s)\}$$
and let $b_0=0$ and $a_{\ell+1}=r$.  Now note that $\textbf{c}\in E_i\Leftrightarrow
\sigma^{\mathbf c}=\sigma_i\Leftrightarrow $ for all $j, \lambda(b_j,b_j+1)$ and
$\lambda(a_j-1,a_j)$ are both in the set
$\{\sigma_i(a_j),\sigma_i(a_j+1),\ldots,\sigma_i(b_j+1)\}\Leftrightarrow $ for all
$j, \lambda(a_j-1,a_j)=\max\lambda(b_{j-1},a_j)$ and
$\lambda(b_j,b_j+1)=\min\lambda(b_j,a_{j+1})$.  Since for any $j$
the chain $x_{b_j+1}<x_{b_j+1}<\ldots<x_{a_{j+1}-1}$ is the unique saturated chain in the
closed interval $[x_{b_j+1},x_{a_{j+1}-1}]$ with increasing labels, the above gives
us:

\begin{fact} \label{FACT} $\textbf{c}\in E_i\Leftrightarrow$ for all $j$ with $0\leq
j\leq \ell, x_{b_j}<x_{b_j+1}<\ldots<x_{a_{j+1}}$ is the unique saturated chain in
the closed interval $[x_{b_j},x_{a_{j+1}}]$ with increasing labels.
\end{fact} 

The following restatement of the previous lemma will be helpful in later sections:  Let $\textbf{c}$ and $\textbf{c}'$ be as in the statement of the previous fact. 
Choose $m\in S$ so that $\lambda(x_{m-1},x_m)>\lambda(x_m,x_{m+1})$, let $\hat{x}_m$ be the element in that interval other than $x_m$ and let $\hat{\textbf{c}}$ be the chain that results from replacing $x_m$ in $\textbf{c}$ with $\hat{x}_m$.  In essence, the previous lemma simply says:  

$$\textbf{c}\in E_i\Rightarrow\hat{\textbf{c}}\in E_i$$

Finally, to prove property (iv), we mimic the construction from section \ref{super}.  Because each $\Sigma_i$ was defined in the same way as in
that case, it suffices to show the following:

\begin{lem}  For any $i<j$ and for any non-maximal chain
$\textbf{c}$ contained in $E_i\cap E_j$, there is some maximal chain $\textbf{d}$ in
$L_j$ containing $\textbf{c}$ as a subchain and such that $\textbf{d}\notin
\Sigma_j$ (i.e., such that $\textbf{d}$ is `old').
\end{lem}

\textbf{Proof:}  Let $\textbf{c}:=x_{\gamma_1}<x_{\gamma_2}<\ldots<x_{\gamma_p}$ be
a non-maximal chain in $E_i\cap E_j$, where the subscripts on the $x$'s indicate
their ranks.  For each $v\in [\ell]$, let $\rho_v=\min\{k:x_k\in \textbf{c}$ and $k\in
\beta_v\}$ and define $\omega_v$ to be the maximum of the same set.  For this proof,
we assume that each $\rho_v$ and $\omega_v$ exists, i.e., that for each interval
$\beta_v$ in $S$ there is some element $x$ of $\textbf{c}$ such that rank$(x)\in
\beta_v$.  The more general case where this is not necessarily true can easily be
worked out with the same technique that we use here.  

Let $\textbf{d}_j:=\hat{0}=y_0<y_1<\ldots<y_r=\hat{1}$, let $\textbf{d}''$ be the chain
obtained by restricting $\textbf{d}_j$ to elements with ranks in
$([r]\cup\{0\})\setminus S$ and adding the elements of $\textbf{c}$.  (Note that
since $\textbf{c}$ is a chain in $L_j$, this process actually gives us back a chain.)  Let
$\textbf{d}'$ be the maximal (in $B_r$) chain obtained by filling in each gap in
$\textbf{d}''$ with the unique maximal chain in that gap with increasing labels, and
let $\textbf{d}$ be the chain $\textbf{d}'$ restricted to elements with ranks in
$S$.  Also, let $\tau\in S_n$ be the label of $\textbf{d}'$.  We wish to show that
$\textbf{d}$ is old, in other words that $\sigma^{\mathbf d}\neq \sigma_j$.  

Let $m=\min\{k\in [r]:\sigma_i(k)\neq \sigma_j(k)\}$.  Then
$\sigma_i(m)<\sigma_j(m)$.  First, suppose that $m$ is not in the interval
$[a_k,b_k+1]$ for any $k$.  Then $m$ is in the interval $[b_{k-1}+2,a_k-1]$ for some
$k$, and so $\sigma_i(m)$ is in the set
$\{\tau(a_k),\tau(a_k+1),\ldots,\tau(\rho_k)\}$.  This implies that $$\tau(a_k)\leq
\sigma_i(m)<\sigma_j(m)\leq \max\lambda(y_{m-1},x_{\rho_k})$$ and so, referring back
to Lemma \ref{lex}, we see that $\sigma_j$ cannot equal
$\sigma^{\mathbf d}$.

To help clarify this argument, the picture below shows an example where $r=9$,
$S=\{3,4,5,6\}$, $\rho_1=4$, $\omega_1=5$, and the chain $\textbf{c}$ is the two
element chain with the solid edge labeled $3$.  In the picture, $\textbf{d}_i$ is
the chain on the left, $\textbf{d}_j$ is the chain on the right, and $\textbf{d}'$
is the labeled chain in the second figure.

\includegraphics[scale=1.1]{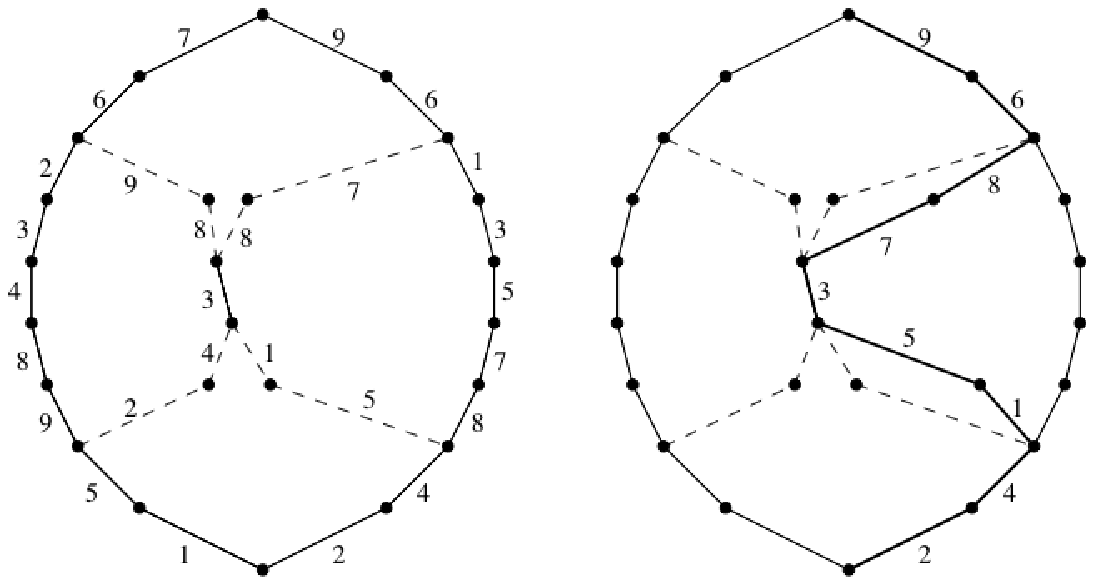}

Now suppose that $m\in [a_k,b_k+1]$ for some $k$.  Because $[a_k,b_k+1]\subseteq
d(\sigma_i)$, it follows that $\sigma_j(m)\notin
\{\sigma_i(a_k),\sigma_i(a_k+1),\ldots,\sigma_i(b_k+1)\}$, and so
$\sigma_j(m)=\tau(n)$ for some $n$ with $\omega_k<n\leq b_k+1$.  Since
$\sigma_i$ is lexicographically earlier than $\sigma_j$, there must be some $q$ such
that $\sigma_i(q)$ is less than every element in the set
$\{\tau(\omega_k+1),\tau(\omega_k+2),\ldots,\tau(b_k+1)\}$.  Either
$\sigma_i(q)\in \{\tau(b_k+2),\tau(b_k+3),\ldots,\tau(a_{k+1}-1)\}$ or $\sigma_i(q)\in
\{\tau(a_{k+1}),\tau(a_{k+1}+1),\ldots,\tau(\rho_{k+1}\}$.  In the first case, $$\tau(b_k+1)\geq\tau(n) >\sigma_i(q)\geq
\min\{\tau(b_k+2),\tau(b_k+3),\ldots,\tau(a_{k+1})\}$$  In the second case,
$$\tau(b_k+1)\geq \tau(n)>\sigma_i(q)\geq \tau(a_{k+1})$$  Either way,
$\tau(n)>\min \{\tau(b_k+2),\tau(b_k+3),\ldots,\tau(a_{k+1})\}$, which means that
$\textbf{d}$ is not a chain in $E_j$, proving property (iv).  $\square$

It is implicit in our work in section 3 that supersolvable lattices are essentially composed of Boolean lattices.  It makes sense, then, that after tackling the rank-selected Boolean lattice, we should find a decomposition for the rank-selected supersolvable lattice.

\section{The Rank-Selected Supersolvable Lattice}

Let $L'$ be a supersolvable lattice of rank $r$ such that $x,y\in L'$ and $x<y$
implies $\mu(x,y)\neq 0$, let $S\subseteq [r-1]$, and let $\lambda$ be a fixed
$S_r$-labeling of $L'$.  We wish to find a convex-ear decomposition of the
rank-selected subposet $L:=L'_S$.  This will generalize our previous two results
(since the Boolean lattice is supersolvable and since $L'_S=L'$ when $S=[r-1]$), but
we will use our previous work to prove this case.

Again, write $S$ as a disjoint union of closed intervals of the integers, such that
the union of no two intervals is an interval: $S=[a_1,b_1]\cup
[a_2,b_2]\cup\ldots\cup [a_\ell,b_\ell]$.  Where appropriate in proofs, let
$\hat{0}=b_0$ and $\hat{1}=a_{\ell+1}$.

As in section \ref{super}, let $\textbf{c}_1,\ldots \textbf{c}_t$ be the chains of $L'$ with decreasing labels under $\lambda$.  For each
$i$ with $1\leq i\leq t$, let $L_i$ be the sublattice of $L'$ generated by
$\textbf{c}_i$ and $\textbf{c}$, where $\textbf{c}$ is the unique maximal chain in
$L'$ with increasing labels (i.e., the M-chain).  As shown in section \ref{super}, each $L_i$
is isomorphic to $B_r$, the Boolean lattice on $r$ elements.  Our approach is
to decompose each $L_i$ based on the results of section \ref{rankbool}, and then combine the
decompositions.  Fix $i$.  Following section \ref{rankbool}, let
$\textbf{c}_i^1,\textbf{c}_i^2,\ldots,\textbf{c}_i^k$ be all maximal chains in $L_i$ (listed in lexicographic order of their labels)
such that $d(\lambda(\textbf{c}_i^j))=S$ for each $j$ (note that $k$ is the same value for each $i$).  Now fix $j$, let
$\textbf{c}_i^j:=\hat{0}=x_0<x_1<\ldots <x_r=\hat{1}$, and define $C_i^j$ as
$\{\textbf{d}:=\hat{0}=y_0<y_1<\ldots <y_r=\hat{1}:n\notin S\Rightarrow x_n=y_n\}$. 
Finally, define the poset $L_i^j$ (not to be confused with the $L_i^j$ mentioned in
section 4) as $\{z\in L: z\in \textbf{c}$ for some $\textbf{d}\in C_i^j\}$ and let
$E_i^j$ be the set of all maximal chains in $L_i^j$ that are not maximal chains in
any $L_i^m$ for $m<j$ or maximal chains in $L_m^\ell$ for any $m<i$ and any $\ell$. 
Note that, often times, $E_i^j$ may be empty.  

For each $i$ and $j$, let $\Sigma_i^j$ be the simplicial complex whose facets are given by the proper maximal chains in $E_i^j$.

\begin{thm}  The sequence of complexes
$\Sigma_1^1,\Sigma_1^2,\ldots,\Sigma_1^k,\Sigma_2^1,\ldots,\Sigma_2^k,\ldots,\Sigma_t^k$,
once we eliminate all $\Sigma_i^j$ such that $\Sigma_i^j=\emptyset$, is a convex-ear
decomposition of $\Delta(L)$.  
\end{thm}

Our work in the previous sections immediately verifies properties i) and ii) of
the decomposition.  

Now let $\textbf{d}':=\hat{0}=y_0<y_1<\ldots<y_r=\hat{1}$ be a maximal chain in
$C_i^j$, and let $\textbf{d}$ be its restriction to elements with ranks in $S$.  Choose
$m\in S$ so that $\lambda(y_{m-1},y_m)>\lambda(y_m,y_{m+1})$, let $\hat{y}_m$ be
the element in the two element interval (in $L_i$) $(y_{m-1},y_{m+1})$ other than
$y_m$, and let $\hat{\textbf{d}}$ be chain resulting from replacing $y_m$ in
$\textbf{d}$ with $\hat{y}_m$. 

\begin{lem} [property (iii)]  Suppose $\textbf{d}\in E_i^j$ for some $i,j$.  Then
$\hat{\textbf{d}}\in E_i^j$.  
\end{lem}

\textbf{Proof:}  Given our work in section \ref{rankbool}, we know that we
cannot have $\hat{\textbf{d}}\in L_i^k$ for some $k<j$.  Suppose that $m\notin
\{a_1,b_1,a_2,b_2,\ldots, b_\ell\}$, and that $\hat{\textbf{d}}$ is a maximal
chain in $L_k^n$ for some $k<i$.  Because of the uniqueness of the maximal chain in the closed interval $[y_{m-1},y_{m+1}]$ with increasing labels (which is the chain
$y_{m-1}<y_m<y_{m+1}$), it must be the case that $\textbf{d}'$ is a maximal chain in
$L_k$, a contradiction.  Now suppose that $m=b_n$ for some $n\in [\ell]\cup\{0\}$. 
Then since $\textbf{d}$ is not in $L_i^k$ for any $k<j$, fact \ref{FACT} tells us that the saturated chain $y_{b_n-1}<y_{b_n}<y_{b_n+1}<\ldots<y_{a_{n+1}+1}$ is the unique saturated
chain in the interval $[y_{b_n-1},y_{a_{n+1}+1}]$ with increasing labels. 
Therefore, if $\hat{\textbf{d}}$ is a maximal chain in $L_k^n$ for some $k<i$,
$\hat{\textbf{d}}'$ (the result of replacing $y_m$ with $\hat{y}_m$ in
$\textbf{d}'$) must be a maximal chain in $L_k$, meaning that $
\hat{\textbf{d}}$ must also be a maximal chain in $L_k$, a contradiction.  The case
when $m=a_n$ for some $n$ is symmetric.  $\square$

As in the previous cases, the following lemma implies property (iv):

\begin{lem} [property (iv)]  Suppose $\textbf{d}'$ is a non-maximal chain that is
both a subchain of some maximal chain in $E_i^j$ for some $i$ and $j$, and a subchain of some
maximal chain in $L_n^p$ where $\langle i, j\rangle $ is lexicographically greater
than $\langle n, p\rangle $.  Let $\textbf{d}$ be the maximal chain in $E_i^j$
obtained by filling in the `gaps' in $\textbf{d}'$ by the unique increasing chain in
those intervals.  Then $\textbf{d}$ is a maximal chain in $L_i^j$ but not $E_i^j$. 
That is, $\textbf{d}$ is an `old' chain in $L_i^j$.
\end{lem}

\textbf{Proof:}  Let $\hat{\textbf{d}}:=\hat{0}=y_0<y_1<\ldots<y_r=\hat{1}$ be the
maximal (in $L'$) chain such that $y_m$ is the element of rank $m$ in $\textbf{d}$
if $m\in S$ and the element of rank $m$ in $\textbf{c}_i^j$ if $m\notin S$.  Suppose
that $\textbf{d}$ is not a maximal chain in $L_i^k$ for any $k<j$.  Then Fact \ref{FACT}
tells us that for every $x,z \in \textbf{d}'$ such that no
$w\in\textbf{d}'$ lies between $x$ and $z$, the chain
$y_{{\rm rank}(x)}<y_{{\rm rank}(x)+1}<\ldots<y_{{\rm rank}(z)}$ is the unique increasing chain in the interval
$[x,z]$.  This guarantees that $\textbf{d}$ is a maximal chain in $L_n$,
meaning that $\textbf{d} \in L_n^k$ for some $k$, proving the lemma.  $\square$

\section{The Rank-Selected Face Poset}\label{faceposet}

Throughout, let $\Sigma$ be a $(d-1)$-dimensional pure shellable simplicial complex.  The result in this section can be seen as motivated by the following proposition, which was proven by Hibi in \cite{hi}.

\begin{prop}  If $\Sigma$ is as above, then the $(d-2)$-skeleton of $\Sigma$ is 2-CM (as defined in Section 8).
\end{prop}

Now let $\Sigma$ have shelling order $F_1,F_2,\ldots, F_t$ and face poset $P_\Sigma$.  For each $i$, let $P_i$ be the face poset of $F_i$ and let $Q_i$ be the face poset of the simplicial complex generated by the first $i$ facets in the shelling order.  

\begin{thm}  Let $S\subseteq [d-1]$.  Then the order complex of the rank-selected poset, $\Delta((P_\Sigma)_S)$, admits a convex-ear decomposition.
\end{thm}

\textbf{Proof:}  Our proof relies on the following fact:  a maximal chain $\textbf{c}$ of $P_i$ is `new' (that is, it is not a maximal chain in $Q_{i-1}$) if and only if some element of $\textbf{c}$ (equivalently, the top element of $\textbf{c}$), when viewed as a face of $\Sigma$, contains the face $r(F_i)$.  Next, we note that each face poset $P_i$ is isomorphic to $B_{d-1}$.  Because of our work in section \ref{rankbool}, we know that $\Delta((P_1)_S)=\Delta((Q_1)_S)$ admits a convex-ear decomposition.  Now assume that $\Delta((Q_{i-1})_S)$ admits a convex-ear decomposition.  We show that we can extend this to a convex-ear decomposition of $\Delta((Q_i)_S)$.  First we define a labeling of $P_i$.  Let $V$ be the set of vertices of $F_i$, and note that any ordering $\phi:V\rightarrow [d-1]$ of $V$ induces an EL-labeling $\lambda$ of $P_i$, since if $y$ covers $x$ then $x\subseteq y$ and $y\setminus x=v\in V$, and we define $\lambda(x,y)=\phi(v)$.  Choose any ordering $\phi$ of $V$ such that if $v\in F_i\setminus r(F_i)$ and $w\in r(F_i)$ then $\phi(v)<\phi(w)$, and let $\lambda$ be the corresponding EL-labeling.  We claim that the convex-ear decomposition of $\Delta((P_i)_S)$ described in section 5 extends the convex-ear decomposition of $\Delta((Q_{i-1})_S)$.  Properties (i) and (ii) are easily verified based on our previous work. 

Taking cues from section \ref{rankbool}, let $\textbf{d}_1,\textbf{d}_2,\ldots,\textbf{d}_\gamma$ be all maximal chains in $P_i$ whose labels have descent set equal to $S$ (where the chains are listed in lexicographic order of their labels).  For each $j$, define $L_j$ and $E_j$ as in section \ref{rankbool}, and let $E_j'$ consists of all maximal chains in $E_j$ that are not maximal chains in $Q_{i-1}$.  We claim that reverse lexicographic order of the chains in $E_j'$ gives a shelling of the associated complex.  To see this, let $\textbf{c}_1$ and $\textbf{c}_2$ be two maximal chains in $E_j'$, where $\lambda(\textbf{c}_1)$ lexicographically precedes $\lambda(\textbf{c}_2)$.  Then $\textbf{c}_2$ comes earlier in the ordering of maximal chains.  Let $x\in F_i$ be the element of highest rank at which the two chains coincide.  If there is an ascent in the label of $\textbf{c}_1$ somewhere it does not coincide with $\textbf{c}_2$ and this ascent is lower than $x$, switching this ascent to a descent as before gives us a chain $\textbf{c}'$ coinciding with $\textbf{c}_1$ above $x$, meaning this chain contains the highest element $y$ of $\textbf{c}_1$.  Since $\textbf{c}_1\in E_j'$, $r(F_i)\subseteq y$ and so $\textbf{c}'\in E_j'$.  So, we only need to worry about the case in which the chains coincide up to $x$.  We consider the case where $d-1\notin S$.  The other case follows easily.  

Let $x<x_1<x_2<\ldots<x_{p}<y$ be the `top' of $\textbf{c}_1$, where rank$(y)=\max S+1$ and $y\in\textbf{d}_j$.  As in section \ref{rankbool} we write $y$ in the expression of $\textbf{c}_1$, even though it is not actually in the chain.  If there is some ascent in the label $\lambda(x,x_p):=\lambda(x,x_1)\lambda(x_1,x_2)\ldots\lambda(x_{p-1},x_p)$, we can switch this to a descent and obtain a chain $\textbf{c}'\in E_j'$ proving our claim.  Otherwise, it must be the case that the label $\lambda(x,x_p)$ is strictly decreasing.  Note that the label $\lambda(x,y)$ cannot be strictly decreasing since $\lambda(\textbf{c}_1)$ is lexicographically earlier than $\textbf{c}_2$.  Therefore, $\lambda(x_{p-1},x_p)<\lambda(x_p,y)$.  Furthermore, since $\textbf{c}_1$ is not a chain in $Q_{i-1}$, it cannot be the case that $\phi^{-1}(\lambda(x_p,y))\in r(F_i)$.  Since $\phi^{-1}(\lambda(x_p,y))\notin r(F_i)$, $\phi^{-1}(\lambda(x_{p-1},x_p))\notin r(F_i)$.  We can therefore switch this ascent to a descent and obtain a maximal chain $\textbf{c}'\in E_j'$.  This proves property (iii).  

The proof of property (iv) carries through just as before, with one other observation needed.  Suppose that $\textbf{c}$ is a non-maximal chain in $Q_{i-1}$ that is a subchain of some maximal chain in $E_j'$, but not a subchain of any maximal chain in $L_m$ for some $m<j$.  We need to produce a maximal chain in $Q_{i-1}$ containing $\textbf{c}$.  But this is immediate since no maximal chain in $Q_{i-1}$ can be in $E_j'$, by definition.  $\square$

The above Theorem cannot be extended to include the case when $d$ is in the subset $S$.  Consider for instance the case when $\Sigma$ consists of a single $2$-dimensional simplex.  Then if $S=\{2,3\}$ it is clear that $\Delta((P_\Sigma)_S)$ has no convex-ear decomposition.  This, however, leads to the following:

\begin{conj}  Suppose that the $(d-1)$-dimensional simplicial complex $\Sigma$ admits a convex-ear decomposition.  Then for any $S\subseteq [d]$, $\Delta((P_\Sigma)_S)$ admits a convex-ear decomposition.
\end{conj}

The following proposition will be needed in Section 8.

\begin{prop}\label{switchlem}  Let $\Sigma$ be as in the statement of the previous theorem, suppose that $S=[d-1]$, let $\textbf{c}:=x_0<x_1<\ldots<x_d$ be a `new' maximal chain in $P_i$ (i.e., $\textbf{c}$ is not a chain in $Q_{i-1}$), and let $\lambda$ be the labeling of the poset $P_i$ as described in the Theorem.  Suppose that $j\in \{1,2,\ldots, d-1\}$ and that $\lambda(x_{j-1},x_j)<\lambda(x_j,x_{j+1})$.  Then if $y$ is the element of $P_i$ in the open interval $(x_{j-1},x_{j+1})$ other than $x_j$ and $\textbf{c}'$ is the chain obtained by replacing $x_j$ in $\textbf{c}$ with $y$, $\textbf{c}'$ is not a chain in $Q_{i-1}$.  In other words, $\textbf{c}'$ is also `new.'
\end{prop}

This proposition is actually true for any $S\subseteq [d-1]$, but we will only use the above case.  

\textbf{Proof:}  If $j<d-1$, then both $\textbf{c}$ and $\textbf{c}'$ contain the element $x_{d-1}$ which, when viewed as a face of $\Sigma$, must contain the face $r(F_i)$ (since $\textbf{c}$ is a `new' chain).  Now suppose that $j=d-1$.  Because $\lambda(y,x_d)=\lambda(x_{d-2},x_{d-1})<\lambda(x_{d-1},x_d)$, we can make the following claim about $y$ and $x_{d-1}$ (when they are viewed as faces of $\Sigma$):  $y$ is obtained from $x_{d-1}$ by removing a vertex $v$ and adding a vertex $v'$, where $\phi(v)<\phi(v')$ (and $\phi$ is as in the proof of the above Theorem).  Since $\textbf{c}$ is new, $r(F_i)\subseteq x_{d-1}$.  Furthermore, since $\phi$ labels all vertices of $r(F_i)$ higher than those vertices not in $r(F_i)$, $r(F_i)\subseteq y$.  Thus, $\textbf{c}'$ is `new.'  $\square$

\section{The Rank-Selected Geometric Lattice}\label{rankgeom}

Let $L$ be a geometric lattice of rank $r$.  In \cite{ns}, Nyman and Swartz showed that
$\Delta(L-\{\hat{0},\hat{1}\})$ admits a convex-ear decomposition.  We briefly describe their technique.

Let $a_1,a_2,\ldots,a_\ell$ be a linear ordering of the atoms of $L$.  The
\textit{minimal labeling} of $L$ labels the edges of the Hasse diagram of $L$ as
follows:  If $y$ covers $x$, then $\lambda(x,y)=\min\{i:x\vee a_i=y\}$.  When
$\textbf{c}$ is a maximal chain, we write $\lambda(\textbf{c})$ to mean its label. 
Bj\"orner has shown (\cite{bj}) that the minimal labeling of $L$ is an EL-labeling.  

In describing the convex-ear decomposition in \cite{ns}, we assume familiarity with some basic notions of matroid theory, including a matroid's lattice of flats and its nbc-bases.  For background on these topics, see for instance \cite{bj}.  

Viewing $L$ as the lattice of flats of a simple matroid $M$, let $B_1,B_2,\ldots,
B_t$ be all its nbc-bases listed in lexicographic order.  For a fixed $i\leq t$,
let $B_i=\{a_{i_1},a_{i_2},\ldots,a_{i_r}\}$ where $a_{i_1}<a_{i_r}<\ldots<a_{i_r}$
under the fixed ordering of the atoms of $L$, and let $b^i_j=a_{i_j}$ for all $j$. 
Fix a permutation $\sigma\in S_r$, and define $\textbf{c}_\sigma^i$ to be the
maximal chain $a_{i_{\sigma(1)}}<a_{i_{\sigma(1)}}\vee a_{i_{\sigma(2)}}<\ldots<
a_{i_{\sigma(1)}}\vee a_{i_{\sigma(2)}}\vee \ldots \vee a_{i_{\sigma(r)}}$.  We
define the basis labeling $\nu_i(\textbf{c})$ of $\textbf{c}$ to be the word
$i_{\sigma(1)}i_{\sigma(2)}\ldots i_{\sigma(r)}$. 

For each $i$ with $1\leq i\leq t$, let $\Delta_i=\{\textbf{c}_\sigma^i:\sigma\in
S_r\}$, and let $E_i$ be all maximal chains in $\Delta_i$ that are not maximal
chains in $\Delta_j$ for any $j<i$.  Finally, let $\Sigma_i$ be the simplicial complex with facets given by the maximal proper chains in $E_i$.

\begin{thm} \cite{ns}  $\Sigma_1,\Sigma_2, \ldots,\Sigma_t$ is a convex-ear decomposition of $\Delta(L)$.  
\end{thm}
     
In proving the above, the authors showed the following (Proposition 4.4):  $$\textbf{c}\in
E_i\Leftrightarrow \nu_i(\textbf{c})=\lambda(\textbf{c})$$

Let $S\subseteq [r-1]$.  We wish to find a convex-ear decomposition of $L_S$. 
The technique will be analogous to the one used in the previous section:

For each $i$ with $1\leq i\leq t$, it is clear that $\Delta_i$ is isomorphic to the
Boolean lattice $B_r$.  Fix $i$, and let
$\textbf{c}_i^1,\textbf{c}_i^2,\ldots,\textbf{c}_i^k$ be all maximal chains in
$\Delta_i$ such that $d(\nu_i(\textbf{c}_i^j))=S$ for all $j$, listed in
lexicographic order of $\nu_i(\textbf{c}_i^j)$.  Now fix $j$, let
$\textbf{c}_i^j:=\hat{0}=x_0<x_1<\ldots<x_r=\hat{1}$, define $C_i^j$ to be the set
$\{\textbf{d}:=\hat{0}=y_0<y_1<\ldots<y_r=\hat{1}: n\notin S\Rightarrow y_n=x_n\}$,
and let $L_i^j=\{z\in L_S:z\in\textbf{d}$ for some $\textbf{d}\in C_i^j\}$. 
Finally, define $E_i^j$ to be the set of all maximal chains in $L_i^j$ that are not
maximal chains in $L_i^n$ for any $n<j$ or maximal chains in $L_n^m$ for any $n<i$,
and let $\Sigma_i^j$ be the associated simplicial complex.  

\begin{thm}  The sequence of complexes
$\Sigma_1^1,\Sigma_1^2,\ldots,\Sigma_1^k,\Sigma_2^1,\ldots,\Sigma_2^k,\ldots,\Sigma_t^k$,
once we eliminate all $\Sigma_i^j$ such that $\Sigma_i^j=\emptyset$, is a convex-ear
decomposition of $\Delta(L_S)$.  
\end{thm}

Note that, as before, properties (i) and (ii) of the decomposition are immediately
verified.  Before verifying the last two properties, however, we need a lemma:

\begin{lem} \label{BLAH} Fix $i$, and suppose
$\hat{\textbf{c}}:=\hat{0}=y_0<y_1<\ldots<y_r=\hat{1}$ is a maximal chain in $E_i$. 
Suppose $\textbf{c}$ is a proper subchain of $\hat{\textbf{c}}$ with the following
property:  anytime $y_p,y_q\in \textbf{c}$ with $y_p<y_q$ and $(y_p,y_q)\cap
\textbf{c}=\emptyset$, $y_p<y_{p+1}<\ldots<y_q$ is the unique maximal chain in
$[y_p,y_q]$ with increasing labels (under $\lambda$).  Then $\textbf{c}$ is not a
subchain of any chain in $\Delta_j$ for $j<i$.  In other words, $\textbf{c}$ is `new.'  
\end{lem}

\textbf{Proof:}  Suppose $\textbf{c}$ is not new, and choose $k$ to be the minimal
$j$ such that $\textbf{c}$ is a subchain of some maximal chain in $\Delta_j$.  Fill
in each gap in $\textbf{c}$ with the unique maximal chain in that interval (in
$\Delta_k$) with an increasing $\nu_k$-label, and call the resulting chain
$\textbf{c}'$.  Assume that $\textbf{c}'\in E_k$.  Then $\nu_k(\textbf{c}')=\lambda(\textbf{c}')$ and by uniqueness of the maximal
chains with increasing $\lambda$-labels in the gaps of of $\textbf{c}$,
$\textbf{c}'=\textbf{c}$, a contradiction since $\textbf{c}\in E_i$ and $k<i$. 
Otherwise, $\textbf{c}'$ is not in $E_k$, meaning it is a maximal chain in
$\Delta_j$ for some $j<i$.  But then $\textbf{c}$ is a subchain of this maximal
chain, contradicting the minimality of $k$.  $\square$
 
Now let $\hat{\textbf{c}}:=\hat{0}=x_0<x_1<\ldots<x_r=\hat{1}$ be a maximal chain in
$\Delta_i$ for some $i$, and let $\textbf{c}$ be its restriction to elements with
ranks in $S$.  Suppose that $m\in S$ with
$\lambda(x_{m-1},x_m)<\lambda(x_m,x_{m+1})$, let $x_m'$ be the element other than
$x_m$ in the two-element open interval $(x_{m-1},x_{m+1})$ in $\Delta_i$, and let
$\textbf{c}'$ be the chain that results from replacing $x_m$ in $\textbf{c}$ with
$x_m'$.

As we noted in the previous section, the following suffices to verify property (iii):

\begin{lem} [property (iii)]  Suppose that $\textbf{c}\in E_i^j$.  Then
$\textbf{c}'\in E_i^j$.
\end{lem}

\textbf{Proof:}  Let $\hat{\textbf{c}}'$ be the chain that results from replacing
$x_m$ with $x_m'$ in $\hat{\textbf{c}}$.  By definition, $\textbf{c}\in
E_i^j\Rightarrow\hat{\textbf{c}}\in E_i$.  The switching lemma from \cite{ns} tells us that $\hat{\textbf{c}}'\in E_i$.  Because each gap in $\textbf{c}$ is filled by
the unique maximal chain in that interval with increasing labels (by Fact
\ref{FACT}), Lemma \ref{BLAH} implies that $\textbf{c}'\in E_i^j$.  $\square$

We also prove property (iv) the same way as in the previous sections:

\begin{lem} [property (iv)]  Fix $i$ and $j$.  Suppose that $\textbf{c}$ is a
non-maximal chain in $L_S$ that is both a subchain of some maximal chain in $E_i^j$
and some maximal chain in $L_n^p$ where $\langle n,p\rangle$ is lexicographically
less than $\langle i, j\rangle$.  Then $\textbf{c}$ is a subchain of some maximal
chain in $L_i^j$ that is not a maximal chain in $E_i^j$.
\end{lem}

\textbf{Proof:}  Again, the algorithm is the same: Fill in the gaps in $\textbf{c}$
with the unique maximal chains in those intervals with increasing $\nu_i$-labels. 
If the resulting chain is new (that is, if it is in $E_i^j$), then $\textbf{c}$ must
be new, by Lemma \ref{BLAH}.  Therefore, the resulting chain cannot be in $E_i^j$. 
$\square$

\section{Flag h-vector Inequalities}

In \cite{ns}, the authors prove the following:

\begin{thm}  (\cite{ns})  Let $L$ be a geometric lattice of rank $r$ with order complex $\Delta=\Delta(L-\{\hat{0},\hat{1}\})$, let $S,T\subseteq [r-1]$, and suppose that $S$ dominates $T$.  Then the flag h-vector of $\Delta$ satisfies:  $h_T\leq h_S$.
\end{thm}

Our goal is to find an analogue of this theorem for face posets of Cohen-Macaulay simplicial complexes.  Although we do not use any of its algebraic properties, we still make use of the Hilbert series of a simplicial complex's Stanley-Reisner ring.  For the definition of this object (as well as further reading on myriad interesting connections between combinatorics and commutative algebra), see \cite{greenstanley}, Section 2.

Let $P$ be a graded poset, and let $\Delta=\Delta(P)$ be its order complex.  Under the fine grading of the face ring $k[\Delta]$, we have 

$$F(k[\Delta],\lambda)=\sum_{F\in \Delta}\prod_{x_i\in F}\frac{\lambda_i}{1-\lambda_i}$$  

We specialize this grading to accommodate the flag h-vector as follows:  identify $\lambda_i$ and $\lambda_j$ whenever the vertices in $\Delta$ to which they correspond have the same rank $r$ (as elements of $P$).  Call this new variable $\nu_r$.  This specialized grading gives us: 

$$F(k[\Delta],\nu)=\sum_{S\subseteq [d]}  f_S\prod_{i\in S}\frac{\nu_i}{1-\nu_i}$$  

\noindent We put this over the common denominator of $\prod_1^d(1-\nu_i)$ to obtain:  

$$F(k[\Delta],\nu)=\sum_{S\subseteq [d]}\frac{f_S \prod_{i\in S}\nu_i\prod_{i\notin S}(1-\nu_i)}{\prod_1^d(1-\nu_i)}=\sum_{S\subseteq [d]}\frac{h_S\prod_{i\in S}\nu_i}{\prod_1^d(1-\nu_i)}$$ 

\noindent Now suppose $\Delta$ triangulates a ball, and let $\Delta'= \Delta-\partial\Delta-\emptyset$.  The following equation is Corollary II.7.2 from \cite{greenstanley}: 

$$(-1)^dF(k[\Delta],1/\lambda)=(-1)^{d-1}\tilde{\chi}(\Delta)+\sum_{F\in \Delta'}\prod_{x_i\in F}\frac{\lambda_i}{1-\lambda_i}$$

\noindent Letting $f_S'$ be the flag f-vector for $\Delta'$, noting that $\tilde{\chi}(\Delta)=0$, plugging in $1/\lambda$ in place of $\lambda$, and specializing to the $\nu$-grading, the previous expression becomes:

$$(-1)^dF(k[\Delta],\nu)=\sum_{S\subseteq [d]}f_S'\prod_{i\in S}\frac{1}{\nu_i-1}$$ 

\noindent Putting the above over the common denominator of $\prod_1^d(\nu_i-1)$ and multiplying by $(-1)^d$ gives us:

$$F(k[\Delta],\nu)=\sum_{S\subseteq [d]}\frac{f_S'\prod_{i\notin S}(\nu_i-1)}{\prod_1^d(1-\nu_i)}$$

\noindent Comparing with our earlier expression for $F(k[\Delta],\nu)$ and noting that the denominators are equal, we have:

$$\sum_{S\subseteq[d]}h_S\prod_{i\in S}\nu_i=\sum_{S\subseteq[d]}f_S'\prod_{i\notin S}(\nu_i-1)$$

\noindent In general, the flag f- and h-vectors satisfy the equation 

$$\sum_{S\subseteq[d]}f_S\prod_{i\notin S}(\nu_i-1)=\sum_{S\subseteq[d]}h_S\prod_{i\notin S}\nu_i$$

\noindent So, we can write the above equation as:

$$
\sum_{S\subseteq[d]}f_S'\prod_{i\notin S}(\nu_i-1)=\sum_{S\subseteq[d]}h_{[d]-S}\prod_{i\notin S}\nu_i
$$

We now apply this above equation to obtain a set of inequalities for the flag h-vector of the face poset of a shellable simplicial complex.

Let $K$ be a $d$-dimensional shellable complex with face poset $P_K$ and shelling order $F_1, F_2, \ldots,  F_t$, and for each $i$ let $P_i$ be the face poset of $F_i$.  Let $A=[d]$, and set $\Sigma = \Delta((P_K)_A)$.  Note that $\Sigma$ is simply the order complex of $P_K$ once we remove the elements corresponding to the facets of $K$ and the element corresponding to the empty set.  Let $\Sigma_1,\Sigma_2,\ldots,\Sigma_t$ be the convex-ear decomposition of $\Sigma$ given in section \ref{faceposet}, and define the sets $E_i'$ as in the proof of Theorem 6.2 (namely, let $E_i'$ consist of all maximal chains of $P_i$ that are not maximal chains in any $P_j$ for $j<i$).  

\begin{lem}\label{switch}  Let $S,T\subseteq [d]$, and suppose that $S$ dominates $T$.  Then for any $i$ there are at least as many maximal chains in $E_i'$ with descent set $S$ (under the labeling described in the proof of Theorem 6.2) as there are with descent set $T$.
\end{lem} 

\textbf{Proof:}  This is an immediate consequence of Lemma \ref{switchlem}.  $\square$

\begin{thm}  \label{flagh}Let $S,T\subseteq[d]$, and suppose that $S$ dominates $T$.  Then the flag h-vector of $\Delta$ satisfies $h_T\leq h_S$.
\end{thm}

\textbf{Proof:}  The argument here is based on the one given in \cite{ns} for geometric lattice order complexes.  First, we note that $h_T(\Sigma_1)\leq h_S(\Sigma_1)$, since the poset associated to $\Sigma$ is just the Boolean lattice $B_d$.  In general, suppose the result holds for $ \Sigma_1\cup\Sigma_2\cup\ldots \Sigma_{k-1}$.  Let $\Omega=\Sigma_1\cup\Sigma_2\cup\ldots\Sigma_{k-1}$, and let $\Sigma_k'=\Sigma_k-\partial\Sigma_k-\emptyset$.  Because $\Sigma_k$ triangulates a ball, we can now use our earlier expression for the flag h-vector of a ball and invoke an argument similar to Chari's in \cite{ch}:

\begin{align*} \sum_{S\subseteq [d]}h_S(\Omega\cup\Sigma_k)\prod_{i\notin S}\nu_i&=\sum_{S\subseteq [d]}f_S(\Omega\cup\Sigma_k)\prod_{i\notin S}(\nu_i-1)\\ 
&=\sum_{S\subseteq [d]}f_S(\Omega)\prod_{i\notin S}(\nu_i-1)+\sum_{S\subseteq [d]}f_S(\Sigma_k')\prod_{i\notin S}(\nu_i-1)\\
&=\sum_{S\subseteq [d]}h_S(\Omega)\prod_{i\notin S}\nu_i+\sum_{S\subseteq [d]}h_{[d]-S}(\Sigma_k)\prod_{i\notin S}\nu_i\\
&=\sum_{S\subseteq [d]}(h_S(\Omega)+h_{[d]-S}(\Sigma_k))\prod_{i\notin S}\nu_i
\end{align*}

Reverse lexicographic order of the maximal chains of $\Sigma_k$ is a shelling, so it follows that $h_S(\Sigma_k)$ is the number of maximal chains of $\Sigma_k$ whose labels have \textit{ascent} set $S$.  Thus $h_{[d]-S}$ counts the number of maximal chains in $P_k$ with descent set $S$.  Since we add at least as many maximal chains whose labels have descent set $S$ as we do maximal chains whose labels have descent set $T$ (Lemma \ref{switch}), the result follows.  $\square$

The previous theorem can now be generalized to face posets of Cohen-Macaulay complexes.

\begin{thm}\label{cmflag}  Let $K$ be a $d$-dimensional Cohen-Macaulay simplicial complex with face poset $P$, and let $\Delta=\Delta(P)$.  Let $S,T\subseteq [d]$, and suppose that $S$ dominates $T$.  Then the flag h-vector of $\Delta$ satisfies $h_T\leq h_S$.
\end{thm}

\textbf{Proof:}  First, note that a linear inequality of the flag h-vector of a complex is equivalent to some linear inequality of the flag f-vector of that complex.  Next we note that, in the case when $\Delta$ is the order complex of the face poset of some complex $K$, a linear inequality of the flag f-vector is equivalent to some linear inequality of the standard $f$-vector of $K$.  To see why this is true, let $S\subseteq [d]$ and write $S$ as a decreasing word:  $a_1,a_2,\ldots, a_m$.  Then $f_S$ is simply the product $b_1b_2\ldots b_m$, where $b_1=f_{a_1}(K)$ and for $i>1$ $b_i=b_{i-1} \binom{a_{i-1}}{a_i}$.  Since Stanley has shown that all linear inequalities involving the f-vector of Cohen-Macaulay complexes are of the form $\sum_ic_ih_i\geq 0$ where each $c_i\geq0$, and since shellable complexes are Cohen-Macaulay, it must be the case that Theorem \ref{flagh} amounts to an inequality of the above form.  Because $h_i(\Delta)\geq 0$ for all $i$ when $\Delta$ is Cohen-Macaulay (see for instance \cite{greenstanley}, pg. 57), it follows that Cohen-Macaulay complexes satisfy the conclusion of Theorem \ref{flagh}.  $\square$ 

We now show that Theorem \ref{cmflag} cannot be extended to include posets whose order complexes are 
Cohen-Macaulay (or 2-CM, for that matter).  We call a graded poset $P$ \textit{Eulerian} if for all $x,y\in P$ 
with $x<y$ we have $\mu(x,y)=(-1)^k$, where $k={\rm{rank}(y)-\rm{rank}(x)}$.  An Eulerian poset whose order complex is Cohen-Macaulay is called \textit{Gorenstein*}.  It can be shown that the order complex of a Gorenstein* poset is 2-Cohen-Macaulay.  For $S\subseteq [n]$, define $w(S)$ to be the set of all $i\in [n]$ such that exactly one of $i$ and $i+1$ is in $S$.  For instance, if $S=\{2,3\}\subseteq [4]$ then $w(S)=\{1,3\}$.  Since Conjecture 2.3 from \cite{eulerian} has been proven by Karu in \cite{ka}, we can rephrase Proposition 2.8 from \cite{eulerian} as:

\begin{prop}  If $S,T\subseteq [n]$ are such that $h_T(\Delta)\leq h_S(\Delta)$ whenever $\Delta$ is the order complex of a Gorenstein* poset then $w(T)\subseteq w(S)$.  
\end{prop}

Now consider $S,T\subseteq [4]$ given by $S=\{1,2\}$ and $T=\{1\}$.  In \cite{ns}, it is shown that $S$ dominates $T$.  However, $w(S)=\{2\}$ and $w(T)=\{1\}$, so $w(T)\nsubseteq w(S)$ and it is clear that we cannot weaken the assumptions of Theorem \ref{cmflag} to include the wider class of Cohen-Macaulay posets (or even 2-CM posets). 

We close this section by mentioning an interesting consequence to Theorem \ref{cmflag}.

\begin{cor}  For each pair of subsets $S,T\subseteq [d-1]$, where $S$ dominates $T$, there exist nonnegative integers $a^{S,T}_1,a^{S,T}_2,\ldots, a^{S,T}_d$ such that, for any $d$-dimensional Cohen-Macaulay $K$ with face poset $P$ and face poset order complex $\Delta=\Delta(P-\emptyset)$,

$$h_S(\Delta)-h_T(\Delta)=\sum_{i=1}^d a^{S,T}_ih_i(K)$$
\end{cor}

\textbf{Proof:}  In the proof of Theorem \ref{cmflag}, we see how linear inequalities of the flag h-vector of $\Delta$ translate to linear inequalities of the h-vector of $K$.  The conclusion of the Theorem tells us that $h_S(\Delta)-h_T(\Delta)\geq 0$ whenever $S$ dominates $T$.  Since all linear inequalities of the h-vector of $K$ must be of the form $\sum_0^da_ih_i(K)\geq 0$, where each $a_i\geq 0$, the corollary follows.  $\square$ 
 
\section{Final Remarks}

As of yet, only a handful of simplicial complexes have been shown to admit convex-ear decompositions.  Other than the ones shown in these paper, these complexes are matroid independence complexes (\cite{ch}) and finite buildings (\cite{sw}).  As shellable simplicial complexes are abundant, the next natural question for such a complex is whether it admits a convex-ear decomposition.

As we already mentioned in section 2, any simplicial complex that admits a convex-ear decomposition is 2-CM.  More intriguing, though, is the following partial converse conjectured by Bj\"orner and Swartz:

\begin{conj}  Let $\Delta$ be a $(d-1)$-dimensional 2-CM simplicial complex.  Then the h-vector of $\Delta$ satisfies $h_i\leq h_{d-i}$ and $h_i\leq h_{i+1}$ for $i<d/2$.  Furthermore, its g-vector is an M-vector.
\end{conj}

It should be noted that although showing a simplicial complex admits a convex-ear decomposition usually amounts to a series of arguments in geometric combinatorics, the underlying machinery of such decompositions takes place in a much more algebraic setting involving the Stanley-Reisner ring of a simplicial complex.  For further reading on such subjects, see \cite{greenstanley} or \cite{st3}.

\textbf{Acknowledgements:}  Authorship of this paper would not have been possible without the inexhaustible guidance and patience of Ed Swartz.  It should also be noted that Vic Reiner initially suggested that order complexes of supersolvable lattices with non-zero M\"obius functions may admit convex-ear decompositions, and that our decomposition of the rank-selected Boolean lattice is based on an unpublished homology basis given by Michelle Wachs.

\end{document}